\input amstex
\documentstyle{amsppt}
\topmatter
\magnification=\magstep1
\pagewidth{5.2 in}
\pageheight{6.7 in}
%\hcorrection{-0.4in}
%\vcorrection{-0.4in}
\abovedisplayskip=10pt
\belowdisplayskip=10pt
\NoBlackBoxes
\title
Sums of powers of consecutive $q$-integers
\endtitle
%Use \endgraf to indicate new paragraph
%\thanks will become a 1st page footnote
\author  Taekyun Kim \endauthor
\affil\rm{{Institute of Science Education,}\\
{ Kongju National University, Kongju 314-701, S. Korea}\\
{e-mail: tkim$\@$kongju.ac.kr}}\\
\endaffil

\abstract{ Let $n, k$ be the positive integers $( k>1 ),$ and let
$S_{n,q}(k)$ be the sums of the $n$-th powers of positive
$q$-integers up to $k-1$: $S_{n,q}(k)=\sum_{l=0}^{k-1}q^l [l]_q^n
$ , where $[l]_q=\frac{q^l-1}{q-1}$. Following an idea due to J.
Bernoulli, we explore a formula for $S_{n,q}(k)$ as follows:
$$ S_{n+1,q}(k)=\frac{1}{n+1}\sum_{i=0}^n
\binom{n+1}{i}\beta_{i,q}q^{ki}[k]_q^{n+1-i}+\frac{(1-q^{(n+1)k})\beta_{n+1,q}}{n+1},$$
where $\beta_{i,q}$ are the $q$-Bernoulli numbers. }
\endabstract
\thanks 2000 Mathematics Subject Classification  11S80, 11B68 \endthanks
\thanks Key words and phrases: $q$-Bernoulli number, $p$-adic $q$-integrals, zeta function, Dirichlet series \endthanks
\rightheadtext{ Taekyun Kim    } \leftheadtext{ Sums of powers of
consecutive $q$-integers  }
%\TagsOnRight
\endtopmatter

\document

\head 1. Introduction \endhead Let $n$ be a natural numbers. There
are formulas such as
$$\aligned
& 1+2+\cdots+n=\frac{n^2+n}{2},\\
&1^2+2^2+\cdots+n^2=\frac{2n^3+3n^2+n}{6}, \\
&1^3+2^3+\cdots+n^3=\frac{n^4+2n^3+n^2}{4}, \cdots.
\endaligned$$
In Korea, these are subjects of the high school mathematics. J.
Bernoulli (1713) first discovered the method which one can produce
those formulae for the sum $\sum_{l=1}^nl^k$ for any natural
numbers. Let $q$ be an indeterminate which can be considered in
complex number field, and for any integer $k$ define the
$q$-integer as $[k]_q=\frac{q^k-1}{q-1}, $ cf. [1,2]. Note that
$\lim_{q\rightarrow 1}[k]_q=k .$ In this paper we evaluate sums of
powers of consecutive $q$-integers.  For any positive integers $n,
k(>1),$ let $S_{n,q}(k)=\sum_{l=0}^{k-1}q^l[l]_q^n.$ Following an
idea due to J. Bernoulli, we explore a formula for $S_{n,q}(k)$ as
follows:
$$ S_{n+1,q}(k)=\frac{1}{n+1}\sum_{i=0}^n
\binom{n+1}{i}\beta_{i,q}q^{ki}[k]_q^{n+1-i}+\frac{(1-q^{(n+1)k})\beta_{n+1,q}}{n+1},$$
where $\beta_{i,q}$ are the $q$-Bernoulli numbers.

\head 2. Sum of the $n$-th powers of positive $q$-integers up to
$k-1$
\endhead

Let $j$ be the positive integers. Then we easily see that
$$[j+1]_q^2-[j]_q^2=([j]_q+q^j)^2-[j]_q^2=q^j(2[j]_q+q^j).\tag 1$$
By (1), we have
$$[k]_q^2=\sum_{j=0}^{k-1}([j+1]_q^2-[j]_q^2)=2\sum_{j=0}^{k-1}q^j[j]_q
+\frac{[2k]_q}{[2]_q}. \tag2$$ Hence,
$$\sum_{j=0}^{k-1}q^j[j]_q=\left([k]_q^2-\frac{[2k]_q}{[2]_q}\right)\frac{1}{2}.\tag3$$
By the same method of Eq.(1), we easily see that
$$[j+1]_q^3-[j]_q^3 =3[j]_q^2q^{j+1}+3[j]_qq^{j} +q^{3j}.\tag4$$
Thus we have
$$\sum_{j=0}^{k-1}[j]_q^2q^{j+1}=\frac{1}{3}[k]_q^3-\frac{1}{2}\left([k]_q^2-\frac{[2k]_q}{[2]_q}\right)-
\frac{1}{3}\frac{[3k]_q}{[3]_q}. \tag5$$ Example. For
$q=\frac{9}{10}$ in (5), we note that
$$\aligned
&(\frac{9}{10})^2+(\frac{9}{10})^3\left(1+\frac{9}{10}\right)^2
+\cdots+(\frac{9}{10})^k\left(1+\frac{9}{10}+\cdots
+(\frac{9}{10})^{k-2}\right)^2+\cdots\\
&=\frac{1}{3}\left(\frac{1}{1-\frac{9}{10}}\right)^3
-\frac{1}{2}\left(\left(\frac{1}{1-\frac{9}{10}}\right)^2-\frac{1}{1-(\frac{9}{10})^2}\right)
-\frac{1}{3}\left(\frac{1}{1-(\frac{9}{10})^3}\right).\endaligned$$
Let $n, k$ be the positive integers $(k>1)$ and let
$S_{n,q}(k)=\sum_{l=0}^{k-1}q^l[l]_q^n. $ Then we note that
$$\sum_{i=0}^{n-1}\binom{n}{i}S_{i,q^{n-i}}(k)=\sum_{l=0}^{k-1}\left([l+1]_q^n-[l]_q^n\right)=[k]_q^n.\tag6$$
By replacing $n$ by $n+1$, we see that
$$[k]_q^{n+1}=\sum_{i=0}^{n}\binom{n+1}{i}S_{i,q^{n+1-i}}(k)
=\sum_{i=0}^{n-1}\binom{n+1}{i}S_{i,q^{n+1-i}}(k)+(n+1)S_{n,q}(k).\tag7$$
Therefore we obtain the following:$$
(n+1)S_{n,q}(k)=[k]_q^{n+1}-\sum_{i=0}^{n-1}\binom{n+1}{i}S_{i,q^{n+1-i}}(k).\tag8$$

\head 3. $q$-analogs of Bernoulli polynomials \endhead
 In this section, we assume $q\in \Bbb C$ with $|q|<1 .$  The
 $q$-Bernoulli polynomials $\beta_{n,q}(x)$ are defined by means
 of the generating function $F_q(t)$ as follows:
$$F_q(t)=e^{\frac{t}{1-q}}\frac{q-1}{\log
q}-t\sum_{n=0}^{\infty}q^{n+x}e^{[n+x]_qt}=\sum_{n=0}^{\infty}\frac{\beta_{n,q}(x)}{n!}t^n,
\text{ $|q|< 1$, $|t|<1$ . }\tag9$$ In the case $x=0,$
$\beta_{n,q}(=\beta_{n,q}(0))$ will be called the $q$-Bernoulli
numbers, cf. [1, 2]. By (9), we easily see that
$$\beta_{n,q}(x)=\sum_{j=0}^{n}\binom{n}{j}q^{jx}\beta_{j,q}[x]_q^{n-j}=\left(\frac{1}{1-q}\right)^n\sum_{k=0}^n\binom{n}{k}
\frac{k}{[k]_q}q^{kx}(-1)^k.\tag10$$ In (9), (10), the
$q$-Bernoulli numbers can be rewritten as
$$\beta_{0,q}=\frac{q-1}{\log q}, \text{ }
(q\beta+1)^k-\beta_{k,q}=\delta_{k,1},\tag11 $$  where
$\delta_{k,1} $ is Kronecker symbol and we use the usual
convention about replacing $\beta^{i}$ by $\beta_{i,q}$, cf.
[1,2]. By (9), (10), and (11), we easily see that
$$-\sum_{l=0}^{\infty}q^{l+k}e^{[l+k]_qt}+\sum_{l=0}^{\infty}q^le^{[l]_qt}
=\sum_{n=1}^{\infty}\left(n\sum_{l=0}^{k-1}q^l[l]^{n-1}\right)\frac{t^{n-1}}{n!}.$$
Thus we have
$$\beta_{n,q}(k)-\beta_{n,q}=n\sum_{l=0}^{k-1}q^l[l]_q^{n-1}.
\tag12$$ By (12), we obtain the following:
$$\sum_{l=0}^{k-1}q^l[l]_q^n=\frac{1}{n+1}\sum_{l=0}^n\binom{n+1}{l}q^{kl}\beta_{l,q}[k]_q^{n+1-l}
+\frac{(1-q^{(n+1)k})}{n+1}\beta_{n+1,q}. \tag13$$ In Eq.(10), we
know that
$$\beta_{n,q}(x)=\sum_{k=0}^n\binom{n}{k}[x]_q^{n-k}\beta_{k,q}\sum_{l=0}^k\binom{k}{l}[x]_q^l(q-1)^l.\tag14$$
By (14), we easily see that
$$\int_{0}^k\beta_{n,q}(x)d[x]_q=\frac{1}{n+1}\left(\beta_{n+1,q}(k)-\beta_{n+1,q}\right)=S_{n,q}(k).
\tag15$$


\Refs \ref\no1 \by T. Apstol \book Introduction to analytic number
theory \bookinfo Undergraduate Texts in Math.\publ
Springer-Verlag, New York\yr 1986\endref \ref \no 2 \by T. Kim
\pages 288-299 \paper $q$-Volkenborn Integration  \yr 2002 \vol 9
\jour Russian J. Math. Phys. \endref

\vskip 0.3cm

\endRefs
\enddocument